\documentclass[10pt]{article}
\usepackage{amsmath,amssymb,amsthm,amscd}
\numberwithin{equation}{section}

\newtheorem{prop}{Proposition}[section]
\newtheorem{theo}[prop]{Theorem}

\newtheorem{rema}[prop]{Remark}

\newtheorem{defi}[prop]{Definition}

\newtheorem{ack}[prop]{Acknowledgment}
\def\begeq{\begin{equation}}
\def\endeq{\end{equation}}

\def\and{\quad{\rm and}\quad}

\def\bl{\bigl(}
\def\br{\bigr)}

\def\<{\langle}
\def\>{\rangle}

\def\lab{\label}

\begin{document}

\title{On Degenerated Monge-Ampere 
Equations over Closed K\"ahler 
Manifolds}

\author{Zhou Zhang
\\Department of Mathematics
\\MIT, Cambridge, MA 02139}
\date{}
\maketitle

\section{Introduction}

The main goal for this note is to prove the following 
theorem which is an improved version of what is stated 
in \cite{t-znote}.  

\begin{theo}

Let $X$ be closed a K\"ahler manifold with (complex) 
dimension $n\geqslant 2$. Suppose we have a holomorphic 
map $F: X\to\mathbb{CP}^N$ with the image $F(X)$ of the 
same dimension. Let $\omega_M$ \footnote{The ``$M$'' is 
the initial letter of ``$model$'' since $\omega_M$ can 
be naturally understood as the model metric of original 
degenerated metric interested in and the degeneration 
information is hiding in the map $F$.} be any K\"ahler 
form over some neighbourhood of $F(X)$ in $\mathbb{CP}^
N$. For the following equation of Monge-Ampere type:  
\begin{equation}
(F^*\omega_M+\sqrt{-1}\partial\bar\partial u)^n=f\Omega,
\end{equation}
where $\Omega$ is a fixed smooth (nondegenerated) volume 
form over $X$ and $f$ is a nonnegative function in $L^p
(X)$ for some $p>1$ with the correct total integral over 
$X$, i.e. $\int_Xf\Omega=\int_X(F^*\omega_M)^n$, \footnote
{This quantity is clearly positive from our assumption.} 
we have the following:  

(1) (Apriori estimate) Suppose $u$ is a weak solution in 
$PSH_{F*\omega_M}(X)\cap L^{\infty}(X)$ of the equation 
with the normalization $\sup_X u=0$, then there is a 
constant $C$ such that $\|u\|_{L^\infty}\le C\|f\|^n_{L^
p}$ where $C$ only depend on $F$, $\omega$ and $p$; 

(2) There would always be a bounded solution for this 
equation;  

(3) If $F$ is locally birational, then any bounded 
solution is actually the unique continuous solution.

\end{theo}

The improvements are in two places: 

i) $X$ closed K\"ahler instead of 
projective; 

ii) in statement $(3)$ about continuity, the 
assumption is weakened a lot.\\ 

It might be worth taking a little time to clarify some 
terminologies appearing in the statement.

First, $u$ is a weak solution means both sides are equal 
as (Borel) measure. The meaning of right hand side is 
classic with $u$ being bounded.

In definition of $L^p(X)$ space, we choose $\Omega$ as the 
volume form. The choice is clearly not so rigid.  

In $(3)$, ``locally birational'' means for a small enough 
neighbourhood $U$ of any point on $F(X)$, each component 
of $F^{-1}(U)$ would be birational to $U$ (under $F$). 
Clearly it would be the case if $F$ is birational itself 
and in fact this is the case with most geometric interests 
as I see it now.\\ 

If $F$ is an embedding, this theorem is proved in \cite{koj98} 
for even more general function $f$. Actually, he only needs 
the $F^*\omega_M$ to represent a K\"ahler class.    

For proving the main theorem above, basically we generalize the 
orginal argument there which makes use of the results in 
\cite{bed-tay} and many other works in pluripotential theory.

\begin{rema}

The discussion in this short note is supposed to be fairly concise. 
We achieve this by taking shortcuts which might make the idea less 
shown. The argument is complete except for freguently referring to 
Kolodziej's \cite{koj98} and \cite{kojnotes}. For greater details 
and more related discussions, we refer the interested readers to 
\cite{zzo}.  

\end{rema}

In this note, plurisubharmonic sometimes means plurisubharmonic 
with respect to some background form. Hopefully, it'll be clear 
from the context.\\ 

\begin{ack}

For projective $X$, the main results of this paper except the 
continuity in (3) of Theorem 1.1 were announced and discussed 
in my previous preprint with Professor Tian. They were also 
presented in a talk by Tian at Imperial College on November 
28, 2005. The general continuity was proved soon in January, 
2006 after a few discussions with Professors S. Kolodziej and 
H. Rossi on approximating plurisubharmonic functions on singular 
spaces. I would like to thank them both for very useful 
discussions. Theorem 1.1 was also presented in my talk at 
Columbia University in February, 2006. A new result in the 
recent preprint by Blocki and Kolodziej allows the current 
generalization to a closed K\"ahler manifold $X$. I really 
appreciate their informing about this result. I would also 
like to thank my advisor, Professor Tian, for bringing up 
his attention to this basic question on complex Monge-Ampere 
equation. 

\end{ack}

\section{Idea of Generalizing Kolodziej's Argument and Preparation}

The (degenerated Monge-Ampere) equation we are considering is the 
following:
$$(\omega_{\infty}+\sqrt{-1}\partial\bar{\partial}u)^{n}=f\Omega$$ 
over $X$, where $\omega_{\infty}=F^{*}\omega_M$ \footnote{This $
\infty$ illustrates the point of view that this degenerated case 
most naturally arises as the limit of nondegenerated case which 
has been used in \cite{t-znote}.}, and $\Omega$ is a volume form 
over $X$. Our goal is to find a bounded (and even continuous) 
solution and get some properties for it, for example, uniqueness. 

\subsection{Idea of Generalization}

Basically, the solution is to be obtained by taking the limit of 
solutions for a family of approximation equations. Of course we 
want the solvability of the approximation equations to be known. 
In our case, the results in \cite{koj98} should be suffcient to 
guarantee this.\\ 

In order to get a limit, we need the apriori ($L^{\infty}$) 
bound for those approximate solutions as in \cite{koj98}. 
In the original argument there, it is the generalized right 
hand side that is treated. We can deal with the right hand 
of the equation above just as well. The main difficulty now 
is of course the degenaracy of $\omega_\infty$ as a K\"ahler 
metric over $X$.

We have a natural family of approximation equations as follows:
$$(\omega_{\infty}+\epsilon \phi+\sqrt{-1}\partial\bar{
\partial}u_{\epsilon})^n=C_{\epsilon}f\Omega$$ 
where $\phi>0$ be closed smooth $(1,1)$-form and $\int_{X}C_{
\epsilon}f\Omega=\int_{X}(\omega_{\infty}+\epsilon\phi)^{n}$ 
for $\epsilon\in (0,1]$. Obviously, we have $C_\epsilon\in 
(1, C]$.

From Kolodziej's result, we have no trouble to find a unique 
continuous solution for each of these equations after 
requiring the normalization $sup_{X}u_{\epsilon}=0$. Now we 
just have to prove that $u_{\epsilon}$'s are uniformly bounded 
from below. 

The difficulty appears since $\omega_{\infty}+\epsilon\phi$ 
is not uniformly positive for $\epsilon\in (0,1]$, i.e., if 
we consider the local potentials in coordinate balls, they will 
no longer be uniformly convex. Thus no matter how good the 
choice is, we do not have a uniform growth of the potential 
when moving away from the origin of the coordinate ball which 
is very crucial for the original argument. 

This one blow seems to completely destroy Kolodziej's argument. 
As we see now, the main reason is that the picture of a 
coordinate (Euclidean) ball in $X$ is a little too local.

The most important observation is that if one chooses a domain 
$V$ which has those degenerated directions of $\omega_\infty$ 
going around inside, we can still have the uniform convex 
local potential, i.e., the values for the very outside part 
are greater than those of the very inside part by a uniform 
positive constant. 

More precisely, if we take a ball in $\mathbb{CP}^{N}$ which 
covers part of $F(X)$, then the preimage of that ball in $X$ 
would be the domain $V$ mentioned above. The potential of $
\omega_{\infty}$ in $V$ be convex in the sense above from the 
positivity of $\omega_M$. Furthermore, we can see the domain 
$V$ is hyperconvex in the usual sense which means we can have 
a continuous exhaustion of the domain, and there are actually 
a lot of nice plurisubharmonic functions over $V$ which can be 
got by pulling back classic functions over the ball in $\mathbb
{CP}^N$.\\ 

There seems to be another problem since the metric $\phi$, 
which is used to perturb the equation, may (should) not have 
a global potential in the domain $V$. But we can deal with 
this by considering plurisubharmonic function in $V$ with 
respect to $\omega_\infty+\epsilon\phi$ for each $\epsilon
\in (0,1]$. In fact we can include the case of $\epsilon=0$ 
in all the argument. The important thing is that our argument 
is uniform for all such $\epsilon$'s.

\begin{rema}

For the apriori estimate in $(1)$ of the main theorem, we only 
need to work with $\omega_\infty$ ($\epsilon=0$). But We'll need 
the estimate uniformly for all $\omega_\infty+\epsilon\phi$ for 
$\epsilon (0, 1]$ for the existence result in $(2)$ of the main 
theorem. The above just says we can treat them together. 

\end{rema} 

We should point out that a global argument will be used below 
which might apparently hide the above idea of considering 
generalized domain $V$. Indeed the punchline is still to study 
the domain $V$. It should be quite natural that after getting 
all the necessary information for the local domains, we can 
patch them up to get for the whole of $X$ just as in 
\cite{kojnotes}.  

\subsection{Preparation} 

Many classic results in pluripotential theory are quite local, 
for example, weakly convergence results, and so can be used 
in our situation automatically. Many definitions also have 
their natural version for the domain $V$ with background 
metric which can't be reduced to potential level globally in 
$V$, for example, relative capacity. Since it is of the most 
importance for us in this work, we give the definition below 
in the case when $V=X$, which is going to be used.

\begin{defi}

Suppose $\omega$ is a (smooth) nonnegative $(1,1)$-form. For any 
(Borel) subset $K$ of $X$, we definite the relative capacity 
of $K$ with respect to $\phi$ as follows: 
$$Cap_{\omega}(K)=sup\{\int_K(\phi+\sqrt{-1}\partial\bar{
\partial}v)^n|v\in PSH_{\phi}(X), ~~-1\leqslant v\leqslant 0
\}.$$ 

\end{defi}

We require $\omega$ to be nonnegative so that $PSH_{\omega}(X)$ 
is not empty. We also point out that usually, it only takes to 
consider compact set $K$ in order to study any set by 
approximation.     

\begin{rema}

Classic results in this business can be found in classic works 
as \cite{le}, \cite{bed-tay}. More recent works as \cite{dem}, 
\cite{kojnotes} might be more convenient as reference.  

\end{rema}

The only thing which is not so trivially adjusted to our situation 
might be comparison principle which is so important and has a 
global feature. There are several ways to deal with this situation.  
One of them seems to be the easiest to describe and in a sense it 
minimizes the modification of Kolodziej's argument for our case. 
So I will use it in this note. The other ways have their own 
interests and will be discussed in \cite{zzo}. 

Now let's state the version of comparison principle we are going to 
use later.

\begin{prop} 
 
For $X$ as above, suppose $u, ~v\in PSH_{\omega}(X)\cap L^{\infty}
(X)$ where $\omega$ is a smooth nonnegative closed $(1,1)$-form, 
then 
$$\int_{\{v<u\}}(\omega+\sqrt{-1}\partial\bar{\partial}u)^n
\leqslant \int_{\{v<u\}}(\omega+\sqrt{-1}\partial\bar{\partial}v)
^n.$$

\end{prop}

This version is slightly different from other more classic versions 
because $X$ may not be projective, $\omega$ may not be positive and 
the functions may not be continuous. The description of justification 
is as follows.\\  

Basically we still just need a decreasing approximation for any 
bounded plurisubharmonic function by smooth plurisubharmonic 
functions according to the argument in \cite{bed-tay}. This is 
not as easy as in Euclidean space where convolution is available. 
And the possible loss of projectivity of $X$ makes it difficult 
to use some other classic results.   

But according to the recent result of Blocki and Kolodziej in 
\cite{blo-koj}, we can have a decreasing smooth approximation 
for plurisubharmonic function over $X$. The approximation result 
need the background form to be positive (i.e., a K\"ahler metric), 
but clearly nonnegative form (as $\omega_\infty$ for us) is 
acceptable when it comes down to comparison principle by simple 
approximation argument.\footnote{We do need $X$ to be K\"ahler.} 
This is also why we can now have $X$ to be just K\"ahler instead 
of projective as stated in \cite{t-znote}.      

\section{Apriori $L^{\infty}$ Estimate}
  
\subsection{Bound Relative Capacity by Measure}

In the following, $\omega$ is a (smooth) nonnegative closed $
(1,1)$-form. Keep in mind that $\omega$ stands for $\omega_
\infty+\epsilon\phi$ for any $\phi\in [0,1]$. The constants do 
not depend on $\epsilon$.  

For $u, ~v\in PSH_\omega(X)\cap L^{\infty}(X)$ with $U(s):=\{u-s
<v\}\neq\varnothing$ for $s\in[S, S+D]$. Also sssume $v$ is valued 
in $[0, C]$. Then $\forall w\in PSH_\omega(X)$ valued in $[-1, 0]$, 
for any $t\geqslant 0$, since $0\leqslant t+Ct+tw-tv\leqslant t+Ct$, 
we have: 
$$U(s)\subset V(s)=\{u-s-t-Ct<tw+(1-t)v\}\subset 
U(s+t+Ct).$$

So we have for $0<t
\leqslant 1$:
\begin{equation}
\begin{split}
\int_{U(s)}(\omega+\sqrt{-1}\partial\bar{\partial}w)^n
&= t^{-n}\int_{U(s)}\bigl( t\omega+\sqrt{-1}\partial\bar{\partial}
(tw)\bigr)^n \\
&\leqslant t^{-n}\int_{U(s)}\bigl( t\omega+\sqrt{-1}\partial\bar{
\partial}(tw)+(1-t)\omega+\sqrt{-1}\partial\bar{\partial}((1-t)v)
\bigr)^n \\
&= t^{-n}\int_{U(s)}\bigl( \omega+\sqrt{-1}\partial\bar{\partial}
(tw+(1-t)v)\bigr)^n \\
&\leqslant t^{-n}\int_{V(s)}\bigl( \omega+\sqrt{-1}\partial\bar{
\partial}(tw+(1-t)v)\bigr)^n \\
&\leqslant t^{-n}\int_{V(s)}\bigl( \omega+\sqrt{-1}\partial\bar{
\partial}(u-s-t-Ct)\bigr)^n \\
&\leqslant t^{-n}\int_{U(s+t+Ct)}(\omega+\sqrt{-1}\partial\bar{
\partial}
u)^n. 
\end{split}
\end{equation}

Comparison principle is applied to get the second to the last 
$\leqslant$. All the other steps are rather trivial from the 
setting. 

Thus we conclude 
$$Cap_\omega(U(s))\leqslant t^{n}\int_{U(s+t+Ct)}
(\omega+\sqrt{-1}\partial\bar{\partial}u)^n.$$ 
from the definition of $Cap_\omega$. 

Let's rewrite this inequality  as: 
$$Cap_\omega(U(s))\leqslant (1+C)^n t^{-n}
\int_{U(s+t)}(\omega+\sqrt{-1}\partial\bar
{\partial}u)^n$$ 
for $t\in (0, min(1, \frac{S+D-s}{1+C})]$. 
Of course, for our purpose, it is always 
safe to assume $\frac{S+D-s}{1+C}<1$.\\ 

Intuitively, the constant $D$ can be seen as 
the gap where the values of $u$ can stretch 
over.  

\subsection{Bound Gap $D$ by Capacity}

We are still in the previous setting. Now 
assume that for any (Borel or compact) 
subset $E$ of $X$, we have: 
$$\int_{E}(\omega+\sqrt{-1}\partial\bar
{\partial}u)^{n}\leqslant A\cdot\frac
{Cap_{\omega}(E)}{Q\bigl( Cap_{\omega}
(E)^{-\frac{1}{n}}\bigr)}$$ 
for some constant $A>0$, where $Q(r)$ 
is an increasing function for positive 
$r$ with positive value. From now on, 
this condition will be denoted by 
Condition $(A)$.\\

The result to be proved in this subsection 
is as follows: 
$$D\leqslant\kappa (Cap_{\omega}(U(S+D)))$$ 

for the following function 
$$\kappa(r)=C_{n}A^{\frac{1}{n}}\bigl(\int_{r^{-\frac{1}{n}
}}^{\infty}y^{-1}(Q(y))^{-\frac{1}{n}}dy+\bigl(Q(r^{-\frac
{1}{n}})\bigr)^{-\frac{1}{n}}\bigr),$$ 
where $C_{n}$ is a positive constant only depending on $n$.\\

The proof is a little technical but quite elementary in spirit. 
We will briefly describe the idea below.\\ 

The previous part gives us an inequality as  
``$Cap \leqslant measure$''.\\ 

Condition $(A)$ gives the other direction 
``$measure \leqslant Cap$''.\\ 

We can then combine them to get some information about the 
length of the interval which comes from $t$ in the inequality 
proved before. The assumption of nonemptiness of the sets is 
needed because we have to divide $Cap_{\omega}(U(\cdot))$ from 
both sides in order to get something for $t$. 

Finally, we can sum all these small $t$'s 
up to control for $D$. \footnote{We use the 
trivial fact that nonemptiness, nonzero 
(Lebesgue) measure and nonzero capacity are 
equivalent for such sets $U(s)$ from the 
fundamental properties of plurisubharmonic 
functions.}\\ 

Of course we'd better use a delicate way to 
carry out all these just in sight of the 
rather complicated final expression of the 
function $\kappa$. It has been done beautifully 
in \cite{koj98}.\\  

Let's point out that in the argument, we do not 
have a positive lower bound for the $t$'s to be 
summed up, so it is important that the inequality 
proved in the previous part holds (uniformly) for 
all small enough $t>0$.

\subsection{Bound Capacity}
 
For $u\in PSH_{\omega}(X)\cap L^{\infty}(X)$ and $u\leqslant 
0$, suppose $K$ is a compact set in $X$ which can well be $X$ 
itself, then there exists a positive constant $C$ such that: 
$$Cap_\omega(K\cap\{u<-j\})\leqslant\frac
{C\|u\|_{L^{1}(V)}+C}{j}.$$ 

\begin{proof}

For any $v\in PSH_\omega(X)$ and valued in $[-1, 0]$, 
consider any compact set $K'\subset K\cap\{u<-j\}$, 
using CLN inequality \footnote{The global version of 
this inequality over $X$ is quite easy to justify in 
sight of the locality of the result.} in \cite{kojnotes}:
\begin{equation}
\begin{split}
\int_{K'}(\omega+\sqrt{-1}\partial\bar{\partial}v)^{n} 
&\leqslant \frac{1}{j}\int_{K}|u|(\omega+\sqrt{-1}
\partial\bar{\partial}v)^{n}\\ 
&\leqslant \frac{C\|u\|_{L^{1}(V)}+C}{j}.\nonumber
\end{split}
\end{equation}

From the definition of relative capacity, this would 
give the inequality above.

\end{proof}  

Now we consider the $L^1$-norm for those approximation 
solutions $u_\epsilon$ (and also the solution $u$ if 
it exists by assumption). The following is just the 
standard Green's function argument. Strictly speaking, 
the computation needs the function to be smooth, but 
we can achieve the final estimate by using approximation 
sequence given by the result in \cite{blo-koj} for our 
situation. So suppose we have the regularity in the 
following.\\ 

For fixed $\epsilon\in [0,1]$, suppose $u_\epsilon(x)=0$ 
and $C>G$ where $G$ is the Green function for the metric 
$\omega_1=\omega_\infty+\phi$. 

Also since $\omega_\infty+\epsilon\phi+\sqrt{-1}\partial
\bar{\partial}u_\epsilon\geqslant 0$, we have 
$$\Delta_{\omega_1}u_\epsilon=\<\omega_1,\sqrt{-1}\partial
\bar{\partial}u_\epsilon\>\geqslant-\<\omega_1, -\omega_
\infty-\epsilon\phi\>\geqslant -C$$
where $C$ is uniform for $\epsilon\in (0,1]$. Basically, 
this tells that there should be no worry for the changing 
background metric. 

Then we have:
\begin{equation}
\begin{split}
0=u_\epsilon(x) 
&= \int_X u_\epsilon{\omega_1}^n+\int_{y\in X}G(x,y)\Delta
_{\omega_1}u_\epsilon\cdot{\omega_1}^n\\
&= \int_X u_\epsilon{\omega_1}^n+\int_{y\in X}\bigl( G(x,y)
-C\bigr)\Delta_{\omega_1}u_\epsilon\cdot{\omega_1}^n\\
&\leqslant \int_X u_\epsilon{\omega_1}^n-C\int_{y\in X}\bigl( 
G(x,y)-C\bigr){\omega_1}^n \\
&\leqslant \int_X u_\epsilon {\omega_1}^n+C.
\end{split}
\end{equation}

This gives the uniform $L^1$ bound for $u_\epsilon$'s by 
noticing they are all nonpositive.\\

Hence we know {\it the set where $u_\epsilon$ has very 
negative value should have (uniformly) small relative 
capacity}.   

\subsection{Conclusion}

Combining all the results above, if we assume  
Condition $(A)$ for some function $Q(r)$ and 
set the function $v$ at the beginning to be 
$0$, we have:
$$D\leqslant\kappa (\frac{C}{D})$$
if $U(s)=\{u<-s\}$ nonempty for $s\in [-2D, -D]$ 
where $C$ is a positive constant.

Furthermore, if we can choose the function $Q(r)$ 
to be $(1+r)^m$ for some $m>0$ so that Condition 
$(A)$ holds, this would imply that the function 
$u$ only take values in a bounded interval since 
$D$ can not be too large \footnote{As $D$ goes to 
$\infty$, $\kappa$ goes to $0$.}. 

That's enough for the lower bound in sight of the 
normalization $sup_X u=0$. The more explicit bound 
claimed in the theorem is not hard to get by 
carefully tracking down the relation. Of course the 
constant $A$ in Condition $(A)$ is fairly involved 
in this business.   
    
\subsection{Condition $(A)$}

In this section, we justify Condition $(A)$ under the 
measure assumption in the main theorem. This part is 
the essential generalization of Kolodziej's original 
argument.\\ 

In our case, $f\in L^{p}$ for some $p>1$, which is the 
measure on the left hand side of Condition $(A)$ from 
the equation we want to solve. For the approximation 
equations, the measures are different, but clearly we 
can bound the $L^p$-norm uniformly. 

Applying H\"older inequality, we know that it suffices 
to prove the following inequality:    
$$\lambda(K)\leqslant A\cdot\bigl(Cap_\omega(K)(1+
{Cap_\omega(K)}^{-\frac{1}{n}})^{-m}\bigr)^q,$$
where $\lambda$ is the smooth measure over $X$ 
and $q$ is some positive constant depending on 
$p$. Obviously, it would be enough to prove:  
$$\lambda(K)\leqslant C_l\cdot{Cap_\omega(K)}^
{l}\cdots\cdots (A)$$ 
for $l$ sufficiently large.\\ 

Of course we have $\lambda(K)<C$, so in fact we 
can get for any nonnegative $l$ if the above is 
true. In the following, we'll consider Condition 
$(A)$ in this form.\\ 

For $\omega$ (uniformly) positive, this can be 
easily reduce to a Euclidean ball. As in 
\cite{kojnotes}, using a classic measure theoretic 
result in \cite{mtsu}, we have: 
$$\lambda(K)\leqslant C\cdot exp\bigl(-\frac{C}
{Cap_\omega(K)^{\frac{1}{n}}}\bigr).$$ 

This is actually stronger than the version above 
after noticing small capacity situation is of the 
main interest.\\    

In the following proof of Condition $(A)$, the essential 
step is to prove the following inequality:
$$\lambda(K)\leqslant C_1\cdot {\epsilon}^{N_1}+C_1
\cdot{\epsilon}^{-N_2}exp\bl\frac{C_2}{{\rm log}
{\epsilon}\cdot Cap_\omega(K)^{\frac{1}{n}}}\br\cdots
\cdots (B),$$ 
for sufficiently small $0<\epsilon<1$. All positive 
constants $C_i$'s do NOT depend on $\epsilon$. This 
$\epsilon$ has nothing to do with the $\epsilon$ 
appearing before in $\omega_\infty+\epsilon\phi$. 

After proving this, by putting $\epsilon={Cap_\omega
(K)}^{\beta}$ for properly chosen $\beta>0$, we can 
justify Condition $(A)$ for any chosen $l$ by noticing 
the dominance of exponential growth over polynomial 
growth.\\

It is easy to notice that we can have uniform 
constants for all $\omega$'s related once we 
get for $\omega_\infty$ from the favorable 
direction of the control we want. And we also 
only need to prove Condition $(A)$ for sets 
close to the subvariety $\{{\omega_\infty}^n=
0\}$ in sight of the results in \cite{kojnotes}.\\     
   
The rest part of this section will be devoted to the 
proof of inequality $(B)$. The following construction 
is of fundamental importance for this goal. 

Let's start with a better description of the map 
$F:X\to F(X)\subset\mathbb{CP}^N$. For simplicity, 
we'll assume here that $F$ provides a birational 
morphism between $X$ and $F(X)$. This assumption 
will be removed at the end. 

Using this assumption, we have subvarieties $Y
\subset X$ and $Z\subset F(X)$ such that $X
\setminus Y$ and $F(X)\setminus Z$ are 
isomorphic under $F$ and $F(Y)=Z$. Clearly $Z$ 
should contain the singular subvariety of $F(X)
$. It's the situation near $Y$ (or $Z$) that is 
of main interest to us.         

Now we use finitely many unit coordinate balls 
on $X$ to cover $Y$. The union of the half-unit 
balls will be called $V$. Then we take two finite 
sets of open subsets depending on $\epsilon>0$ 
as follows:\\  

$\{U_i\}$, $\{V_i\}$, with $i\in I$, finite 
coverings of $V\setminus W$, where $W$ is the 
intersection of $\epsilon$-neighbourhood of 
$Y$ \footnote{That's a neighbourhood of $Y$ 
correspondent to the intersection of balls 
of radius $\epsilon$ in $\mathbb{CP}^N$ 
covering $Z$.} with $F(X)$, such that each 
pair $V_i\subset U_i$ is in one of the chosen 
unit coordinate balls. Moreover, $F(U_i)$ and 
$F(V_i)$ are the intersections of $F(X)$ with 
balls of sizes $\frac{1}{2}{\epsilon}^C$ and 
$\frac{1}{6}{\epsilon}^C$ where $C>0$ are 
chosen to be big enough to justify the above 
construction.  

Clearly $|I|$ is controlled by $C\cdot{\epsilon}
^{-N_2}$.\\     

For any compact set $K$ in $V$, we have the 
following computation:
\begin{equation}
\begin{split}
\lambda(K) 
&\leqslant \lambda(W)+\sum_{i\in I}{\lambda(K\cap\bar{V}_i)}\\
&\leqslant C\cdot{\epsilon}^{N_1}+\sum_{i\in I}{C\cdot exp{
\bl -\frac{C}{{Cap(K\cap\bar{V}_i, U_i)}^{\frac{1}{n}}}\br}}\\
&\leqslant C\cdot{\epsilon}^{N_1}+\sum_{i\in I}{C\cdot exp{
\bl\frac{C}{{\rm log}{\epsilon}\cdot{Cap_{\omega_\infty}(K\cap
\bar{V}_i)}^{\frac{1}{n}}}\br}}\\
&\leqslant C\cdot{\epsilon}^{N_1}+\sum_{i\in I}{C\cdot exp{\bl
\frac{C}{{\rm log}{\epsilon}\cdot{Cap_{\omega_\infty}(K)}^{
\frac{1}{n}}}\br}}\\
&\leqslant C\cdot{\epsilon}^{N_1}+C {\epsilon}^{-N_2}\cdot exp
{\bl\frac{C}{{\rm log}{\epsilon}\cdot{Cap_{\omega_\infty}(K)}^
{\frac{1}{n}}}\br}.\nonumber
\end{split}
\end{equation} 

That's just what we want. $C_1$ and $C_2$ are used 
in the original statement of $(B)$ since the $C$'s 
at different places have different affects on the 
magnitude of the final expression. Of course, the 
same $C$ for each term in the big sum have to be 
really the same constant. In the following, we 
justify the computation above. The only nontrivial 
steps are the second and third ones.\\ 

The second one is the direct application of $(\star)$, 
the classic result in $\mathbb{C}^n$ as $V_i$ and $U_i$ 
are in one of the finitely many unit coordinate balls 
which clearly can be taken as the unit Euclidean ball 
in a uniform way.\\ 

The third step uses the following inequality: 
$$Cap(K\cap\bar{V_i}, U_i)\leqslant C\cdot 
(-{\rm log}{\epsilon})^n \cdot Cap_{\omega_\infty}
(K\cap\bar{V_i}).$$     

This result also has its primitive version in classic 
pluripotential theory for domains in $\mathbb{C}^n$.   
Extension of plurisubharmonic function is all what we 
need to prove it as described below.  

For any $v\in PSH(U_i)$ valued in $[-1, 0]$. If we 
can ``extend'' this function to an element $-C{\rm 
log}\epsilon\cdot\tilde{v}$ where $\tilde{v}$ is 
plurisubahrmonic with respect to $\omega_\infty$ 
valued in $[-1,0]$ over $X$, and also make sure 
that the measures $(\sqrt{-1}\partial\bar{\partial}
v)^n$ and $(\omega_\infty+\sqrt{-1}\partial\bar{
\partial}\tilde{v})^n$ are the same over $\bar{V_i}$, 
then this would clearly imply the inequality above 
from the definition of relative capacity.\\      

The construction will be done mostly on $F(X)$. The 
function $v$ can be considered over $F(U_i)$. We'll 
``extend'' it to a neighbourhood $F(X)\setminus O_i$ 
in $\mathbb{CP}^N$ where $O_i$ is a neighbourhood of 
$\bar{V_i}$ in $U_i$. 

Let's first extend it locally in $\mathbb{CP}^N$. We 
can safely assume that the construction happens in 
(finite) half-unit Euclidean balls in $\mathbb{CP}
^N$ which cover the variety $Z$ and have $\omega_M$ 
defined on the correspondent unit balls. $\omega_M$ 
are can be expressed in the level of potential and 
so the construction is merely about functions. 

Consider the plurisubharmonic function function 
$$h=\bigl({\rm log}(\frac{36|z|^2}{{\epsilon}^{2C}})
\big))^{+}-2,$$ 
where the upper $+$ means taking maximum with $0$, 
on the unit ball in $\mathbb{CP}^N$ but with the 
coordinate system $z$ centered at the center of $F
(V_i)$. It's easy to see the pullback of this 
function, still denoted by $h$, is plurisubharmonic 
and $max(h, v)$ on $U_i$ is equal to $v$ near $
\bar{V_i}$ and equal to $h$ near $\partial{U_i}$. 
So this function extends $v$ to the preimage of 
the unit ball in $\mathbb{CP}^N$ while keeping 
the values near $\bar{V_i}$. 

Now we want to extend further to the whole of $X$. 
We still work on $F(X)\subset\mathbb{CP}^N$. And 
it only left to extend the function $h$ for the 
remaining part where the value is less restrictive. 

$|h|$ is bounded by $-C\cdot{\rm log}\epsilon$ in 
the unit ball. So we can have 
$$\sqrt{-1}\partial\bar{\partial}h=-C\cdot{\rm log}
\epsilon\cdot (\omega_M+\sqrt{-1}\partial\bar{
\partial}H)$$ 
for $H$ plurisubharmonic with respect to $\omega_M$ 
valued in $[-1, 0]$ in the unit ball. Then using 
the same argument as in \cite{kojnotes}, we can 
extend $H$ to (uniformly bounded) $\tilde{H}\in PSH
_{\omega_M}(O)$, where $O$ is a neighbourhood of $F
(X)$, using the positivity of $\omega_M$. Finally we 
just take $\tilde{v}=F^*\tilde{H}$.\\ 

This ends the argument for the case when $F:X\to F
(X)$ is a birational map.\\ 

Now we want to remove the birationality condition. In 
fact, after removing proper subvarieties $Y$ and $Z=F
(Y)$ of $X$ and $F(X)$ respectively, we can have $F:X
\setminus Y\to F(X)\setminus Z$ is a finitely-sheeted 
covering map, since the map is clearly of full rank 
there and the finiteness of sheets can be seen by 
realizing the preimage of any point in $F(X\setminus 
Z)$ is a finite set of points.  

Then it's easy to see the argument before would still 
work in this situation. Basically, we still have the 
construction before. Now the only difference is that 
now the numbers of small pieces $U_i$ and $V_i$ need 
to be multiplied by the number of sheets and this 
won't affect the argument too much.

Hence we get the apriori $L^\infty$ bound in general.     

\section{Existence of Bounded Solution}

Now we discuss the existence of a bounded solution. 
As suggested in Section 2, approximation of the 
orginal equation in the main theorem is used. The 
(uniform) apriori estimate got in the previous 
section would give us the uniform $L^\infty$ bound 
for the approximation solutions $u_\epsilon$ there. 

Thus exactly the same argument of taking limit as 
used in \cite{koj98} can be applied for our case to 
get a bounded solution for the original equation. 

\begin{rema}

Indeed, $u_\epsilon$ is essentially decreasing as 
$\epsilon\to 0$ which will make the limit easier to 
take. More details about this and some uniqueness 
results for bounded solutions will be provided in 
\cite{zzo}.     

\end{rema}

\section{Continuity of Bounded Solution and 
Stablity Result}

In fact, we can prove that a bounded solution for 
the original equation is actually continuous with 
just a little more assumption.\\  

The argument is almost in the same line as the 
argument in \cite{koj98} (using the $L^\infty$ 
argument above). 

We still just have to analyze the situation around 
a carefully chosen point. But now the point might 
be in the set $\{{\omega_\infty}^n=0\}$. So in 
order to have convexity of the local potential of 
$\omega_\infty$, we have to consider the domain 
in $X$ which is the preimage of a ball in $\mathbb
{CP}^n$ under the map $F$. In other words, we do 
have to consider the domain $V$ mentioned earlier.  

For such a domain, where $\omega_\infty+\sqrt{-1}
\partial\bar{\partial}u=\sqrt{-1}\partial\bar{
\partial}U$ for a bounded plurisubharmonic 
function $U$, we don't have convolution which 
gives a decreasing smooth approximation of any 
plurisubharmonic function. 

That's where we'll use the additional assumption of 
$F$ being locally birational and the classic result 
in \cite{for-nar}. The reason for requiring this 
local birationality would be clear below. 

Basically, we want to be able to push forward the 
function $U$ above to the singular domain in $F(X)$ 
in the most straightforward way (without averaging, 
etc.). The blowing-down picture would work. But 
it'll be OK if there are several conponents in $X$ 
correspondent to the same singular domain in $F(X)$ 
as we only need to treat each of the components. In 
a sense, we just do not want any branching. Notice 
that the most interesting case of $F$ being 
birational to its image falls right into this 
category.   

Then it's quite straightforward to see the pushforward 
function $F_*U$ is a weak plurisubharmonic function. 
One needs to prove that for any holomorphic map from 
unit disk in $\mathbb{C}$ to $F(X)$, the pullback of 
$F_*U$ is subharmonic. The idea is to reduce to map 
holomorphic with image in subvariety of $F(X)$ with 
decreasing dimension. 
 
Now the classic result in \cite{for-nar}, pointed out 
to me by Professor Kolodziej, tells us that we can 
(locally) extend $F_*U$ to a plurisubharmonci function 
in a ball of $\mathbb{CP}^N$ which would be enough for 
us to go through local argument (for a properly chosen 
$V$). Basically now we can again use convolution to get 
the desired approximation.\\ 

Thus we can justify continuity for bounded solutions as 
in \cite{koj98}. The punchline is as follows. Suppose $
\{U_j\}$ is the sequence of smooth plurisubarmonic 
functions constructed above which are defined on a 
neighbourhood slightly larger than $V$ which decreases 
to $U$ pointwisely. Then by the construction in 
\cite{koj98}, which is very local and can be easily 
adjusted to our case, we can prove the sets $\{U+c<U_j
\}$ are nonempty and relatively compact \footnote{For 
the relative compactness of the sets, strictly speaking, 
we have to use another function which is constructed 
from $U$ linearly instead $U$ itself in the definition 
as in \cite{koj98}. It's a little too tedious to describe 
the details here.} inside $V$ for all $c\in (0, a)$ for 
$a>0$ and $j>j_0$.

The argument for $L^{\infty}$ estimate before gives 
$$\frac{a}{2}\leqslant\kappa(Cap(\{U+\frac{a}{2}<U_j\}, 
V)).$$  

\begin{rema}

Let's point out that we have to use local argument to get 
this inequality. For this, the justification of comparison 
principle would be different from the global case since 
there is now boundary to consider. Also the locality of 
the approximation above using the result in \cite{for-nar} 
seems to make it difficult to get the approximations for 
two functions for the same domain.  

We can deal with this problem using another way mentioned 
in the discussion about comparison principle. Basically, 
we use a different definition of relative capacity where 
only continuous functions are considered in taking the 
supremum. Then together with the approximation we have for 
$U$, we can go through the local argument for the sets $
\{U+c\leqslant U_j\}$. At this moment, it might be a little 
too distracting to carry out all the details which will 
appear in \cite{zzo}.   

\end{rema}

We also notice that the relative capacity of the set $\{U+
\frac{a}{2}<U_j\}$ would go to $0$ as $j\to\infty$. This is 
can be justified by the decreasing convergence and $\frac{a}
{2}>0$. At last, we can draw the contradiction by letting $j$ 
goes to $\infty$ in the equality above because the right hand 
side is going to $0$. 

\begin{rema}

The continuity can be achieved without the additional 
assumption if we can have a decreasing approximation 
of the solution by functions in $PSH_{\omega_\infty}
(X)\cap C^\infty(X)$. There are all kinds of results 
in this direction, but somehow I feel semi-positivity 
of $\omega_\infty$ won't be sufficient.   

\end{rema}

Finally, let's point out that the stability result for 
bounded (hence continuous) solutions in \cite{kojnotes} 
also hold in the current case. The original proof there 
works for us without any essential change.\\ 

We would like to mention that the stability argument 
there can almost be directly used to consider merely 
bounded solution. This is not pointless here since the 
continuity result needs a little bit more assumption 
than boundedness result as for now. But there is an 
inequality used there which bounds the measures of 
mixed terms from two plurisubharmonic functions by the 
measures from each of them. This inequality seems to 
be hard to prove for merely bounded functions. More 
discussions for this can also be found in \cite{zzo}. 

\section{Application}  

The most useful application of the results above would 
be the $L^\infty$ estimate. Combining it with other 
estimates from $PDE$ methods (for example, maximum 
principle), we can further understand the solution. 

Let's illustrate this by the application in 
K\"ahler-Ricci flow. This has been discussed 
in \cite{t-znote} which focuses on the maximum 
principle argument. 
 
Consider the following K\"ahler-Ricci flow
\begin{equation}
\lab{eq:kahlerricciflow}\frac{\partial{\tilde
{\omega}_t}}{\partial t}= -{\rm Ric}(\tilde
{\omega}_t)-\tilde{\omega}_t,~~~~\tilde{\omega}
_0=\omega_{0},
\end{equation}
where $\omega_0$ is any given K\"ahler metric 
and ${\rm Ric}(\omega)$ denotes the Ricci form 
of $\omega$, i.e., in the complex coordinates, 
${\rm Ric}(\omega)=\sqrt{-1}R_{i\bar{j}}dz^i
\wedge d\bar z^{j}$ where $(R_{i\bar{j}})$ is 
the Ricci tensor of $\omega$. Let $\omega_
\infty=-{\rm Ric}\Omega$ for a volume form 
$\Omega$. Set $\omega_t=\omega_\infty+e^{-t}
(\omega_0-\omega_\infty)$ and $\tilde{\omega}
_t=\omega_t+\sqrt{-1}\partial\bar{\partial}u$, 
we can put (\ref{eq:kahlerricciflow}) on the 
level of potential and more in the Monge-Ampere 
setting as:
$$(\omega_t+\sqrt{-1}\partial\bar{\partial}u)^n
=e^{\frac{\partial u}{\partial t}+u}\Omega.$$
We have seen in \cite{t-znote} that the 
right-handed side has a uniform $L^p$-bound for 
all $t$ with any $1<p\leqslant\infty$.

Now suppose $[\omega_\infty]=-K_X$ is nef and 
big. Hence it would be semi-ample and fall right 
into the setting of the main theorem. 

Though generally in $\omega_t=\omega_\infty+e^{-t}
(\omega_0-\omega_\infty)$, it may not be true that 
$\omega_0-\omega_\infty>0$. Combining with the 
degenerated lower bound of $u$, we can still have 
uniform $L^\infty$ estimate for $u(t,\cdot)$ with 
$t\in [0, \infty)$ simply by using part of $\omega_
\infty$ in the front to dominate the second term. 
\footnote{In fact, the uniqueness result in \cite
{t-znote} allows us to only consider the case when 
$\omega_0>\omega_\infty$.} This would give us the 
boundedness of the limit and the continuity follows 
as well since we can choose the map $F$ to be 
birational to its image for this case.\\ 

There would be generalization and more application 
in \cite{zzo}.

\end{document}